\documentclass[12pt,a4paper]{article}
\usepackage[cp1251]{inputenc}
\usepackage[english,russian]{babel}
\usepackage{graphicx}
\usepackage{amssymb}
\usepackage{amsmath}
\usepackage{amsthm}
\usepackage{euscript}

\newtheorem{theorem}{Теорема}

\newtheorem{rem}{Замечание}

\newenvironment{remark}{\begin{rem}\rm}{\end{rem}}

\newtheorem{ex}{Пример}

\newcommand{\RP}[0]{\mathbb{RP}}
\newcommand{\CP}[0]{\mathbb{CP}}
\newcommand{\CC}[0]{\mathbb{C}}
\newcommand{\R}[0]{\mathbb{R}}

\newcommand{\be}[0]{\begin{equation}}
\newcommand{\ee}[0]{\end{equation}}
\newcommand{\bez}[0]{\begin{equation*}}
\newcommand{\eez}[0]{\end{equation*}}

\newcommand{\ep}[0]{$\hspace{\fill} \square$}

\begin{document}

\title{ Пример прямой без вещественных точек,
лежащей в пересечении трех
комплексных квадрик}
\author{И.Н. Шнурников}
\date{}

\maketitle

\paragraph{Постановка задачи.}
В шестимерном евклидовом пространстве $\R^6$ с координатами $x=(x_1, \dots, x_6)$  множество $Q^3$ задано как совместная поверхность уровней функций
$$
Q^3=\{f_1(x)=d_1, f_2(x)=d_2, f_3(x)=d_3\} \subset \R^6
$$
для функций
\begin{gather*}
f_1(x)=x^2_1+x^2_2+x^2_3+x^2_4+x^2_5+x^2_6\\
f_2(x)=x_1x_4+x_2x_5+x_3x_6\\
f_3(x)=c_1x^2_1+c_2x^2_2+c_3x^2_3+c_4x^2_4+c_5x^2_5+c_6x^2_6,
\end{gather*}
где $c_i$ и $d_j$ --- некоторые вещественные числа. Необходимо

\begin{itemize}
  \item Найти все значения параметров $c_i$ и $d_j$, при которых множество $Q^3$ является гладким многообразием (подмногообразием $\R^6$).
  \item Найти многообразие $Q^3$ с точностью до гомеоморфизма.
\end{itemize}

\paragraph{Схема решения.} А.Б.~Жеглов предложил использовать теоремы Коллара \cite{Kollar} для того, чтобы предъявить некий список многообразий, которому обязано принадлежать $Q^3$. А  именно:

 (а) Многообразие $Q^3$ является множеством вещественных точек пересечения трех квадрик в комплексном проективном пространстве
$$
Q^3_c=\bar{Q}_1\cap \bar{Q}_2\cap \bar{Q}_3 \subset
\CP^6,
$$
где квадрики $\bar{Q}_i$ задаются однородными многочленами от семи комплексных переменных, соответствующими функциям $f_i$. Для этого достаточно проверить, что все вещественные точки многообразия $Q^3_c$ аффинны.

(б) Находятся условия на параметры $c_i,\quad d_j$, при которых пересечение комплексных квадрик будет алгебраическим подмногообразием в $\CP^6$ ({\it  полным пересечением}).

(в) Из известной формулы
 для размерности множества одномерных линейных подпространств
 (см. \cite[ch. 13,
\S. 6, th. 1]{Chodg}) следует, что в $Q^3_c$ есть однопараметрическое
 семейство комплексных прямых.
Определим проекцию
$$
p: Q^3_c \to \CP^2
$$
зафиксировав прямую $m$, лежащую в $Q^3_c$,
рассмотрев $\CP^2$ как множество квадрик, содержащих
$Q^3_c$, и точке $x$ поставив в соответствие квадрику, содержащую
плоскость $(x,m)$.

(г) Отображение $p$ не определено на прямой $m$ и на пересекающих ее прямых. Зато отображение $p$ можно доопределить до отображения $\tilde p$ раздутия $\tilde
Q^3_c$ многообразия $Q^3_c$ вдоль $m$ и пересекающих ее прямых.

(д) Теорема Коллара для отображения
$$
\tilde p: \tilde Q^3_c \to
\CP^2
$$
утверждает, что вещественные точки многообразия $\tilde Q^3_c$ есть одно из следующих многообразий: $\RP^3,S^3,S^2\times S^1$,
  их связные суммы,
 многообразия Зейферта с не более 6 особыми слоями, линзы с $p,q \leq 6$.

(е) Если прямая проекции $m$ и все пересекающие ее прямые не имеют вещественных точек, то раздутие
 не касается вещественных точек и поэтому множества вещественных точек (точек с вещественными однородными координатами) многообразий $\tilde Q^3_c$ и $Q^3_c$ совпадают.

 В итоге, если доказать существование прямой $m$ на $Q^3_c$, такой что на ней и на пересекающих ее прямых, лежащих на $Q^3_c$, нет вещественных точек, то тогда многообразие $Q^3$ будет гомеоморфно одному из указанному в пункте (д).

\paragraph{Результаты работы и дальнейшие перспективы.} В данной работе найдены условия на параметры квадрик, при которых
множество $Q^3$ является подмногообразием $\R^6$ и в $\CC^6$.  Предъявлены прямая $l \subset Q^3_c$ и несколько условий, при выполнении которых на прямой $l$ нет вещественных точек. Для дальнейшего решения задачи необходимо
\begin{itemize}
  \item выяснить, при каких значениях параметров квадрик выполняются указанные условия,
  \item проверить, что прямые, пересекающие $l$ и лежащие на $Q^3_c$, не содержат вещественных точек.
\end{itemize}

\paragraph{История возникновения и актуальность задачи.} Задача возникла в теории динамических систем при описании системы, аналогичной системе Эйлера движения трехмерного тела с закрепленной точкой, и поэтому названной "движением четырехмерного твердого тела"$, $ см. \cite{Bolsinov_99}. Рассматривается гамильтонова система на четырехмерном симплектическом многообразии $M^4$ ---  поверхности уровня функций $f_1(x)=d_1, f_2(x)=d_2$ с гамильтонианом $f_3(x)$. Известно \cite{Bolsinov_99}, что если параметры
удовлетворяют соотношению
$$
c_1c_4(c_2+c_5-c_3-c_6)+c_2c_5(c_3+c_6-c_1-c_4)+c_3c_6(c_1+c_4-c_2-c_5)=0
,$$ то гамильтонова система интегрируема по Лиувиллю и тем самым {\it изоэнергетическая поверхность} $Q^3$ расслаивается на торы, а слоение описывается с помощью молекул и их меток, см. обзоры работ А.~Т.~Фоменко и др. в \cite{Bolsinov_99}. В известных до сей поры примерах интегрируемых систем топология изоэнергетических поверхностей относительно простая, поэтому было бы интересно найти, чему диффеоморфно многообразие $Q^3$ в случае неинтегрируемости.

\paragraph{Невырожденность в $\R^6$ и в $\CC^6$.}

\begin{theorem} Множество уровня трех вещественных функций
$$
Q^3_R=
\{ f_1=d_1,\quad f_2=d_2,\quad f_3=d_3 \}\in\R^6
$$
является
трехмерным вещественным многообразием тогда и только тогда, когда

(а) $d_1>2|d_2|$  и

(б) не существует такого вещественного числа $b$, которое
удовлетворяло бы всем трем уравнениям:
\begin{equation*}
\begin{split}
  &     1) \quad (c_1-a)(c_4-a)=b^2,\\
   &   2) \quad (c_2-a)(c_5-a)=b^2,\\
   &   3) \quad (c_3-a)(c_6-a)=b^2,  \\
\end{split}
\end{equation*}
и хотя бы одному неравенству:
\begin{equation*}
\begin{split}
& 1) \quad\frac {c_1-a}{bd_2}\geq \frac{1+(\frac {c_1-a}{b})^2}{d_1}, \\
& 2) \quad\frac {c_2-a}{bd_2}\geq \frac{1+(\frac {c_2-a}{b})^2}{d_1}, \\
& 3) \quad\frac {c_3-a}{bd_2}\geq \frac{1+(\frac {c_3-a}{b})^2}{d_1}, \\
\end{split}
\end{equation*}
 где $a=\frac {d_3-2bd_2}{d_1}.$
\end{theorem}

\proof
 Эти условия получаются как условия, при которых матрица Якоби производных
 функций $f_1, \quad f_2, \quad f_3$ по переменным $X_1,\dots, X_6$ размера
  $3\times 6$ невырождена во всех точках множества $Q^3_R.$

Обратно, если условие (а) не выполняется, то в матрице Якоби первые
2 строчки зависимы или множество $Q^3_R$ --- пустое.

Если условие (б) не выполняется, то ранг матрицы Якоби равен 2 в
 точке с координатами $x_1,\dots,x_6,$ где

$$
x_4=(\frac {c_4-a}{b})x_1,\quad x_5=(\frac {c_4-a}{b})x_2,\quad x_6=(\frac
{c_4-a}{b})x_3,
$$

 а числа $x_1, x_2, x_3$ являются решением системы:
\begin{equation*}
\begin{split}
 & d_1=x_1^2\left(1+\frac {(c_1-a)^2}{b^2}\right)+x_2^2\left(1+\frac {(c_2-a)^2}{b^2}\right)+x_3^2\left(1+\frac {(c_3-a)^2}{b^2}\right),\\
 & d_2=x_1^2\left(\frac {(c_1-a)}{b}\right)+x_2^2\left(\frac {(c_2-a)}{b}\right)+x_3^2\left(\frac {(c_3-a)}{b}\right), \\
 & d_3=x_1^2\left(c_1+c_4\frac {(c_1-a)^2}{b^2}\right)+x_2^2\left(c_2+c_5\frac {(c_2-a)^2}{b^2}\right)+x_3^2\left(c_3+c_6\frac {(c_3-a)^2}{b^2}\right). \\
   \end{split}
\end{equation*}
\ep

\begin{theorem}
 Множество уровня трех комплексных функций
$$
Q^3_C=\{ f_1=d_1,\quad f_2=d_2,\quad f_3=d_3 \}\in\CC^6
$$
является
трехмерным комплексным многообразием тогда и только тогда, когда

(а) $d_1\neq 2|d_2|$  и

(б) не существует такого комплексного числа $b$, которое
удовлетворяло бы всем трем уравнениям:
\begin{equation*}
\begin{split}
  &     1) \quad (c_1-a)(c_4-a)=b^2,\\
   &   2) \quad (c_2-a)(c_5-a)=b^2,\\
   &   3) \quad (c_3-a)(c_6-a)=b^2,  \\
\end{split}
\end{equation*}
\end{theorem}

\begin{remark}
 Невырожденность в $\CP^6$ получается
расмотрением 7 аффинных карт, в каждой из которых условия
невырожденности будут аналогичны условиям невырожденности в
$\CC^6.$
\end{remark}

\paragraph{Прямая без вещественных точек.} На многообразии  $Q^3_c$ будем искать прямую в виде
$$
x_i=a_i+tb_i\quad \text{для}\quad i=1,2,\dots 6,
$$
где коэффициенты $a_i, b_i$ --- это комплексные числа, а $t$ ---
 это комплексная переменная. Наложим условия на коэффициенты прямой:
$$
a_4=\lambda a_1,\quad a_5=\lambda a_2,\quad a_6=\lambda a_3\quad\text{и}\quad b_4=\mu
b_1,\quad b_5=\mu b_2,\quad b_6=\mu b_3.
$$
После подстановки параметрического задания прямой в три уравнения,
  задающие многообразие $Q^3_c,$ получим девять уравнений (при $t^2, t$
  и свободном члене для каждого из трех уравнения $f_j=d_j$),
   но из-за дополнительной симметрии, возникающей из наложенных
    условий на коэффициенты прямой, независимых
уравнений будет семь, переменных же восемь:
$a_1,a_2,a_3,b_1,b_2,b_3,\lambda,\mu$; поэтому положим число
$\lambda$ равным корню уравнения $x+\frac 1x=\frac{d_1}{d_2},$ и
число $\mu$ равным корню уравнения

\begin{gather*} \sqrt{c_2-c_3+\mu^2(c_5-c_6)}\sqrt{c_3-c_1+\mu^2(c_6-c_4)}(d_3-\frac{d_2}{\lambda}(c_1+\lambda^2c_4))(c_2-c_3+\lambda^2(c_5-c_6))+\\
 + (c_1-c_3+\lambda^2(c_4-c_6))\left(\frac{d_2}{\lambda}(c_2+\lambda^2c_5)-d_3\right)(c_1-c_2+\lambda \mu(c_4-c_5))+ \\
+ \sqrt{c_2-c_3+\mu^2(c_5-c_6)}\sqrt{c_1-c_2+\mu^2(c_4-c_5)}\left(\frac{d_2}{\lambda}(c_2+\lambda^2c_5)-d_3\right)(c_1-c_3+\lambda \mu(c_4-c_6))+\\
 +\sqrt{c_3-c_1+\mu^2(c_6-c_4)}\sqrt{c_1-c_2+\mu^2(c_4-c_5)}
(d_3-\frac{d_2}{\lambda}(c_1+\lambda^2c_4))(c_2-c_3+\lambda \mu(c_5-c_6))=\\
 = 0.
\end{gather*}
Теперь находятся числа
\begin{equation*}
\begin{split}
&b_1=\sqrt{c_2-c_3+\mu^2(c_5-c_6)},\\
&b_2=\sqrt{c_3-c_1+\mu^2(c_6-c_4)},\\
&b_3=\sqrt{c_1-c_2+\mu^2(c_4-c_5)},
\end{split}
\end{equation*}
и находятся числа $a_1,\quad a_2,\quad a_3$.

\begin{theorem}
 Предположим, что
 \begin{itemize}
   \item $
\frac{d_2}{\lambda}(c_2+\lambda^2c_5)-d_3 \neq 0,
$
   \item не выполняется
одновременно пара условий $c_1=c_2$ и $c_4=c_5,$
   \item не выполняется
одновременно пара условий $c_1=c_3$ и $c_4=c_6,$
   \item у уравнения на
$\mu$ есть корень, отличный от $\lambda$ и от $-\lambda,$
\item не все
числа $a_1,a_2,a_3$ вещественны.
 \end{itemize}
Тогда для любого комплексного числа $t$ шесть чисел $a_j+tb_j$ не
могут быть одновременно вещественными.
\end{theorem}

\proof
 Предположим противное, то есть что числа
$a_i+tb_i$ и $\lambda a_i+t\mu b_i$ одновременно вещественны, тогда
три числа $\lambda (a_i+tb_i)$ вещественны. Вычтем вторые числа из
третьих (для всех $i=1,2,3$), получим вещественность чисел
 $t(\lambda-\mu)b_i.$ Деля полученные числа одно на другое (они не равны нулю по сделанным предположениям), получим вещественность чисел $\frac{b_2}{b_1}$ и $\frac{b_3}{b_1}.$
Числа $b_i$ удовлетворяют уравнению $b_1^2+b_2^2+b_3^2=0.$ Из
вещественности чисел $\frac{b_2}{b_1}$ и $\frac{b_3}{b_1}$ следует,
что $b_1=b_2=b_3=0.$ Однако подстановка $b_1=b_2=b_3=0$ в уравнение
на $\mu$ даст уравнение

$$
(c_1-c_3+\lambda^2(c_4-c_6))\left(\frac{d_2}{\lambda}(c_2+\lambda^2c_5)-d_3\right)(c_1-c_2+\lambda \mu(c_4-c_5))=0
$$
Каждый из трех множителей не может быть равен нулю по сделанным
предположениям.
\ep

\bigskip
\begin{remark} Можно заметить, что для почти всех значений параметров
на $Q^3_c$ существует вещественно одномерное множество вещественных
точек, сквозь которые проходят прямые $m\in Q^3_c$. Множество же
самих прямых вещественно двумерно, однако прямые могут густо
пересекаться, поэтому применение теоремы
 Коллара нетривиально, и будет осуществлено с помощью найденной явно прямой.
\end{remark}

\medskip
{\it Благодарен А.Б. Жеглову за ценные указания, ссылки и внимание к работе.}
 
Работа выполнена при поддержке гранта  Правительства
РФ по постановлению N 220, договор No. 11.G34.31.0053.

\begin{flushleft} ЯрГУ, лаб. дискретной и вычислительной
геометрии им. Делоне.\\
E-mail: shnurnikov@yandex.ru
\end{flushleft}

\end{document}